# ESTIMACIÓN Y ANALISIS DE SENSIBILIDAD PARA EL COEFICIENTE DE DIFUSIVIDAD EN UN PROBLEMA DE CONDUCIÓN DE CALOR

Guillermo Federico Umbricht[&], Diana Rubio[*]

Centro de Matemática Aplicada, ECyT Universidad Nacional de San Martin

[&]guilleungs@yahoo.com.ar   [*]diarubio@gmail.com

**Resumen**

En este artículo se estudia la estimación del coeficiente de difusividad de una barra metálica homogénea conociendo valores temporales de temperaturas en un punto intermedio. Para ello, se analiza el estado transitorio de un problema de conducción de calor en un hilo conductor totalmente aislado de longitud $l$, con condiciones de contorno constantes. El problema puede ser modelado con una ecuación en derivadas parciales parabólica, imponiendo condiciones de borde de tipo Dirichlet.

Para el análisis del error, consideramos la temperatura simulada en un punto de la barra para distintos instantes y estimamos el coeficiente de difusividad utilizando técnicas usuales de problemas inversos. Las ecuaciones se discretizaron mediante diferencias finitas centradas y usando funciones de MATLAB.

Incluimos un estudio analítico y numérico de la sensibilidad de la temperatura con respecto al coeficiente de difusividad. Los experimentos numéricos muestran una muy buena precisión en las estimaciones.

**Palabras Clave**: **Ecuación de calor, Condiciones de Dirichlet, Estimación de parámetros, Sensibilidad.**

# ESTIMATION AND SENSITIVITY ANALYSIS OF THE DIFFUSIVITY COEFFICIENT FOR A PROBLEM OF HEAT CONDUCTION

## Abstract

The aim of this article is to discuss the estimation of the diffusivity coefficient of a homogeneous metal rod from temperature values at a fixed point in the bar for different time instants. The time-dependent problem of heat conduction is analyzed in an insulated conductor wire of length l considering constant boundary conditions. The problem is modeled by a parabolic partial differential equation, imposing Dirichlet boundary conditions. We consider simulated temperature values at a point of the bar for different time instants and estimate the coefficient of diffusivity using usual techniques for solving inverse problems. For the discretization of the equation we consider a finite difference centered scheme. We include an analytical and numerical study of the sensitivity of the temperature function with respect to the coefficient of diffusivity. Numerical experiments show very good accuracy in the estimates.

**Key words:** Heat equation, Dirichlet conditions, Parameter estimation, Sensitivity.

## Elección del Tema

Los problemas inversos son cada vez más estudiados en diferentes disciplinas y su creciente interés se puede observar en los numerosos trabajos que aparecen en la bibliografía. Uno de los motivos de este interés creciente se relaciona directamente con la necesidad/ importancia de predecir situaciones futuras o estimar parámetros de un determinado sistema físico a partir de observaciones actuales.

## Definición del Problema

Estudiamos el problema de estimar numéricamente el coeficiente de difusividad, $\alpha^2 [m^2/s]$, en la ecuación de calor de una dimensión con condiciones de borde de tipo Dirichlet, a partir de un conjunto de $K$ mediciones $u(x_0; t_1), \ldots, u(x_0; t_K)$ en un punto medio, $x_0$, de la barra.

El problema puede describirse mediante la ecuación en derivadas parciales de tipo parabólica

$$\begin{cases} u_t(x;t) = \alpha^2 \cdot u_{xx}(x;t) & 0 < x < l, \quad t > 0 \\ u(x;0) = \Phi(x) & 0 < x < l \\ u(0;t) = k_1 & t > 0 \\ u(l;t) = k_2 & t > 0 \end{cases} \quad (1)$$

donde

$u(x;t)$ representa la temperatura en la posición x en el instante t (ºC),

$\Phi(x)$ representa el perfil inicial de temperaturas (ºC),

$k_1$ representa la fuente caliente (ºC),

$k_2$ representa la fuente fria (ºC),

$l$ representa la longitud de la barra $(m)$.

Esquemáticamente:

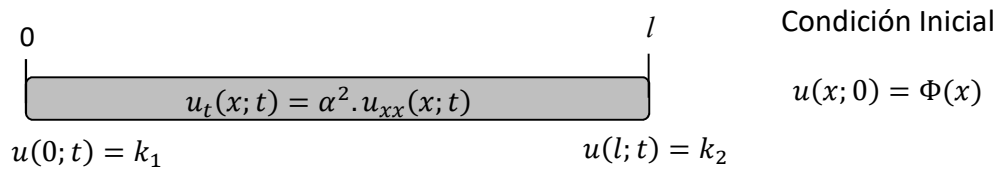

## Antecedentes y Justificación del Estudio

El problema de identificación de parámetros ha sido muy estudiado y ha recibido considerable atención de muchas investigaciones actuales debido esencialmente a que se han encontrado disimiles y múltiples aplicaciones en los distintos campos de la ciencia tales como la conducción de calor [10], la identificación de fisuras [11], la teoría electromagnética [3], la prospección geofísica [2], la detección de contaminantes [6], la tomografía y la tomografía de impedancia eléctrica [1] y determinación de tumores en un tejido biológico [9], entre otros.

La ecuación de calor, o de difusión, se utiliza para describir numerosos procesos de transferencia de masa y energía, es muy estudiada y analizada por áreas como la termodinámica y la termoquímica [4].

En la literatura se pueden encontrar una gran cantidad de algoritmos numéricos para la identificación de distintos parámetros en ecuaciones en derivadas parciales, particularmente en la ecuación de calor. Esto indica la necesidad constante de encontrar nuevas técnicas de estimación para este problema.

## Limitaciones y Alcances del trabajo

En este trabajo nos concentramos en un problema de identificación de un parámetro de la ecuación de conducción de calor con condiciones de borde de tipo Dirichlet (ver [7]).

## Objetivos

- Estudiar analítica y numéricamente la sensibilidad de la solución del problema con respecto al coeficiente de difusividad.
- Estimar numéricamente el coeficiente de difusividad a partir de datos en un punto intermedio de la barra.

## Hipótesis

Estudiar analítica y numéricamente la sensibilidad de la solución con respecto al coeficiente de difusividad y estimar el valor de dicho coeficiente mediante técnicas usuales para la resolución de problemas inversos. Esto nos permitirá identificar el material desconocido con el cual está construido el hilo conductor.

## Desarrollo del trabajo

### 1- Análisis Teórico

#### i. Solución analítica

En el siguiente lema se obtiene la expresión analítica de la solución del problema de conducción de calor en 1 variable espacial, considerando condiciones de borde de tipo Dirichlet.

**Lema 1**

Consideremos el problema parabólico

$$\begin{cases} u_t(x;t) = \alpha^2 \cdot u_{xx}(x;t) & 0 < x < l, \quad t > 0 \\ u(x;0) = \Phi(x) & 0 < x < l \\ u(0;t) = k_1 & t > 0 \\ u(l;t) = k_2 & t > 0 \end{cases} \quad (1)$$

donde $u(x;t)$ es la temperatura de la barra de longitud $l$ en la posición $x$ en el instante $t$, suponiendo constantes los parámetros físicos involucrados: $(\alpha^2, k_1, k_2)$. La solución al problema se puede expresar como:

$$u(x;t) = k_1\left(1 - \frac{x}{l}\right) + k_2\frac{x}{l} + \sum_{n=1}^{\infty} A_n \cdot e^{\frac{-n^2\pi^2}{l^2}\alpha^2 t} \cdot sen\left(\frac{n\pi}{l}x\right) \quad (2)$$

donde

$$A_n = \frac{2}{l}\int_0^l \left(\Phi(x) - k_1\left(1 - \frac{x}{l}\right) - k_2\frac{x}{l}\right) sen\left(\frac{n\pi}{l}x\right) dx \quad (3)$$

**Demostración:**

Consideramos primero el caso estacionario, es decir,

$$u_t^{EST}(x,t) = 0$$

ya que la temperatura depende sólo de la posición $x$ dentro de la barra, no varía con el tiempo. Para esto basta considerar un instante $t$ suficientemente grande. Luego, de la ecuación (1) resulta $u_{xx}^{EST}(x;t) = 0$ o equivalentemente,

$$u^{EST}(x,t) = Ax + B,$$

y usando las condiciones de borde se obtiene:

$$u^{EST}(x,t) = k_1\left(1 - \frac{x}{l}\right) + k_2\frac{x}{l}.$$

El problema (1) se puede resolver agregando a $u^{EST}(x,t)$ una función que tenga dependencia temporal, es decir,

$$u(x;t) = u^{EST}(x,t) + \varphi(x;t)$$

y por lo tanto:

$$u(x;t) = k_1\left(1 - \frac{x}{l}\right) + k_2\frac{x}{l} + \varphi(x;t). \tag{4}$$

A partir de la ecuación y las condiciones de borde que debe satisfacer $u(x;t)$ se obtiene la siguiente ecuación con condiciones inicial y de borde para $\varphi(x;t)$:

$$\begin{cases} \varphi_t(x;t) = \alpha^2 \cdot \varphi_{xx}(x;t) & t > 0; \quad 0 < x < l \\ \varphi(0;t) = 0 & t > 0 \\ \varphi(l;t) = 0 & t > 0 \\ \varphi(x;0) = \Phi(x) - k_1\left(1 - \frac{x}{l}\right) - k_2\frac{x}{l} & 0 < x < l \end{cases} \tag{5}$$

Utilizamos el método de separación de variables, es decir, suponemos que existen funciones $X, T$ tal que:

$$\varphi(x;t) = X(x)T(t) \qquad 0 < x < l, \quad t > 0$$

y reemplazando en la ecuación de calor, obtenemos:

$$X(x) \cdot T'(t) = \alpha^2 \cdot X''(x) \cdot T(t) \implies \frac{X''(x)}{X(x)} = \frac{T'(t)}{\alpha^2 T(t)} = \xi; \qquad \text{con } \xi \text{ constante.}$$

Se tiene entonces un problema de autovalores de la forma $X''(x) - \xi X(x) = 0$ en la variable espacial que tiene solución no trivial sólo en el caso $\xi < 0$. Imponiendo las condiciones de borde se obtiene:

$$X_n(x) = sen\left(\frac{n\pi}{l}x\right) \quad \text{con } n \in \mathbb{Z}.$$

Luego resulta $T'(t) - \xi \alpha^2 T(t) = 0$, de donde

$$T_n(t) = e^{-\frac{n^2\pi^2}{l^2}\alpha^2 t}, \quad n \in \mathbb{Z}.$$

Usando el principio de superposición tenemos que

$$\varphi(x;t) = \sum_{n=1}^{\infty} A_n \cdot e^{-\frac{n^2\pi^2}{l^2}\alpha^2 t} \cdot sen\left(\frac{n\pi}{l}x\right) \tag{6}$$

donde $A_n$ son los coeficientes de Fourier y se definen como

$$A_n = \frac{2}{l}\int_0^l \varphi(x;0) sen\left(\frac{n\pi}{l}x\right) dx.$$

De las condiciones para $\varphi(x;t)$ dadas en (5) se tiene la ecuación (3), esto es,

$$A_n = \frac{2}{l}\int_0^l \left(\Phi(x) - k_1\left(1 - \frac{x}{l}\right) - k_2\frac{x}{l}\right) sen\left(\frac{n\pi}{l}x\right) dx.$$

Por lo tanto, reemplazando (6) en (4) se tiene (2):

$$u(x;t) = k_1\left(1 - \frac{x}{l}\right) + k_2\frac{x}{l} + \sum_{n=1}^{\infty} A_n\, sen\left(\frac{n\pi}{l}x\right) \cdot e^{-\frac{n^2\pi^2}{l^2}\alpha^2 t}.$$

∎

### ii. Análisis de sensibilidad de la temperatura con respecto a $\alpha^2$

Una de las herramientas utilizadas para medir la influencia de un parámetro en los resultados de un determinado modelo es la función de sensibilidad. Algunos ejemplos concretos de utilización y análisis teórico pueden verse en [5],[7],[8].

En esta subsección, calculamos la sensibilidad de la solución con respecto al coeficiente de difusividad ($\alpha^2$). Con este análisis pretendemos mostrar que pequeños errores en el cálculo de $\alpha^2$ no provocarán grandes errores en el cálculo de la temperatura.

Denotamos
$$S(x;t) = \frac{\partial u(x;t)}{\partial \alpha^2}$$
la función de sensibilidad de $u$ con respecto a $\alpha^2$, que cuantificará la dependencia del parámetro en la solución. De aquí que la función $u$ se puede representar como $u = u(x, t, \alpha^2)$.

Para simplificar la notación, en el análisis que sigue se omiten las dependencias de $x$ y $t$ y se define
$$c = \alpha^2.$$
Utilizando el desarrollo de Taylor de primer orden alrededor de un valor determinado $c_0$ obtenemos:
$$u(c) = u(c_0) + u_c(c_0).(c - c_0) + \frac{u_{cc}(\zeta)}{2}.(c - c_0)^2$$
donde el subíndice $c$ denota la derivada parcial con respecto a $c$ y $\zeta$ es un valor intermedio entre $c$ y $c_0$. Por el teorema del resto de Taylor, el término $\frac{u_{cc}(\zeta)}{2}.(c-c_0)^2$ es el error cometido al aproximar la función $u(c)$ por su desarrollo en polinomio de Taylor de primer orden alrededor de $c_0$.

Ya que $S = \frac{\partial u}{\partial c}$ se tiene
$$u(c) = u(c_0) + S(c_0).(c - c_0) + \frac{S_c(\zeta)}{2}.(c - c_0)^2$$
Suponiendo que para $\zeta$ entre $c$ y $c_0$ la función $S_c$ está acotada y que $|c - c_0|$ es pequeño, podemos escribir
$$u(c) \cong u(c_0) + S(c_0).(c - c_0)$$
o, equivalentemente
$$|u(c) - u(c_0)| \cong |S(c_0)|.|c - c_0|.$$
Bajo estas condiciones, podríamos concluir que perturbaciones en el valor de $c_0$ se amplificarán o disminuirán en $u$ de acuerdo al factor $|S(c_0)|$.

Para realizar un estudio numérico de la sensibilidad derivamos con respecto a $\alpha^2$ las ecuaciones que satisface la solución $u(x;t)$ y consideramos sólo las dependencias espacial y temporal para la función de sensibilidad. Obtenemos así las ecuaciones que debe satisfacer la función de sensibilidad
$$\begin{cases} S_t(x;t) = \alpha^2.S_{xx}(x;t) + u_{xx}(x;t) & 0 < x < l, \quad t > 0 \\ S(x;0) = 0 & 0 < x < l \\ S(0;t) = 0 & t > 0 \\ S(l;t) = 0 & t > 0 \end{cases}$$

**Lema 2**

La solución del problema parabólico

$$\begin{cases} S_t(x;t) = \alpha^2 \cdot S_{xx}(x;t) + u_{xx}(x;t) & 0 < x < l, \quad t > 0 \\ S(x;0) = 0 & 0 < x < l \\ S(0;t) = 0 & t > 0 \\ S(l;t) = 0 & t > 0 \end{cases} \quad (6)$$

donde $S(x;t)$ es sensibilidad de la temperatura con respecto a $\alpha^2$ está dada por:

$$S(x;t) = -\frac{\pi^2}{l^2} t \sum_{n=1}^{\infty} n^2 A_n \, sen\left(\frac{n\pi}{l} x\right) \cdot e^{-\frac{n^2 \pi^2}{l^2} \alpha^2 t} \quad (7)$$

donde $A_n$ son los mismos coeficientes de la función temperatura $u(x;t)$ dados en (3).

**Demostración:**

Sigue de manera inmediata que (7) satisfacen las condiciones inicial y de frontera dadas en (6):

$$S(x;0) = 0$$
$$S(0;t) = 0$$
$$S(l;t) = 0$$

Resta ver que se satisface la ecuación diferencial

$$S_t(x;t) = \alpha^2 \cdot S_{xx}(x;t) + u_{xx}(x;t)$$

Derivando la expresión (7) para $S$ y (2) para $u$ se tiene

$$S_t(x;t) = -\frac{\pi^2}{l^2} \sum_{n=1}^{\infty} n^2 A_n \, sen\left(\frac{n\pi}{l} x\right) \cdot e^{-\frac{n^2 \pi^2}{l^2} \alpha^2 t} \left[-\frac{n^2 \pi^2 \alpha^2}{l^2} t + 1\right]$$

$$S_{xx}(x;t) = \frac{\pi^2}{l^2} t \sum_{n=1}^{\infty} n^2 A_n \, sen\left(\frac{n\pi}{l} x\right) \cdot e^{-\frac{n^2 \pi^2}{l^2} \alpha^2 t} \left(\frac{n\pi}{l}\right)^2$$

$$u_{xx}(x;t) = -\sum_{n=1}^{\infty} A_n \, sen\left(\frac{n\pi}{l} x\right) \cdot e^{-\frac{n^2 \pi^2}{l^2} \alpha^2 t} \left(\frac{n\pi}{l}\right)^2$$

De donde se deduce de manera inmediata que:

$$S_t(x;t) - \alpha^2 \cdot S_{xx}(x;t) - u_{xx}(x;t) = 0$$

∎

### iii. Estimación del coeficiente de difusividad $\alpha^2$

La estimación de parámetros es básicamente un problema de optimización. Se trata de encontrar el valor que mejor ajusta la solución basada en un conjunto de datos o mediciones de la solución.

Formalmente, se define el funcional

$$J(\alpha^2) = \sum_{k=1}^{K} |u(x_0, t_k; \alpha^2) + \varepsilon_k - u_k|^2$$

donde $u_k$, $k = 1, \ldots K$ representan los datos o mediciones de la temperatura en la posición $x_0$ en los instantes $t_1, \ldots, t_k, \ldots, t_K$ y $u(x_0, t_k; \alpha^2) + \varepsilon_k$ $k = 1, \ldots K$ representan los valores simulados (mediciones) tomando $\alpha^2$ como coeficiente de difusión siendo $\varepsilon_k$ los errores de medición.

Se supone que existe un valor $\alpha_0^2$ que corresponde al valor real, es decir, que satisface $u(x_0, t_k; \alpha^2) = u_k$ y es el valor que se quiere estimar. Además, los errores de medición se suponen independientes, idénticamente distribuidos, con media 0 y varianza $\sigma^2$.

Bajo estas suposiciones, el estimador de para $\alpha^2$ se puede obtener como

$$\widehat{\alpha^2} = \operatorname*{argmin}_{\alpha^2} J(\alpha^2) = \operatorname*{argmin}_{\alpha^2} \sum_{k=1}^{K} |u(x_0, t_k; \alpha^2) + \varepsilon_k - u_k|^2$$

### 2- Desarrollo numérico

Con la finalidad de estudiar numéricamente las soluciones de la ecuación de calor y analizar la estabilidad de dichas soluciones con respecto a coeficiente de difusividad, discretizamos la ecuación utilizando diferencias finitas centradas en $x$ y a derecha en $t$.

Definimos la partición
$$\mathcal{P} = \{(x_i, t_j), \quad i = 1, \ldots, N, j = 1, \ldots, M \quad x_i \in \mathcal{P}_x, t_j \in \mathcal{P}_t\}$$
donde
$$\mathcal{P}_x = \{0 = x_0 < x_1 < \cdots < x_i < \cdots < x_N = l, \quad x_i = (i-1)\Delta x, \quad i = 1, \ldots, N\}$$
$$\mathcal{P}_t = \{0 = t_0 < t_1 < \cdots < t_j < \cdots < t_M, \quad t_j = (j-1)\Delta t, \quad j = 1, \ldots, M\}.$$

Las ecuaciones discretizadas con las respectivas condiciones iniciales y de borde son:

$$\begin{cases} u(x_i, t_{j+1}) = \varepsilon. [u(x_{i+1}, t_j) + u(x_{i-1}, t_j)] + (1 - 2\varepsilon) u(x_i, t_j), \\ u(x_i; 0) = \Phi(x) & i = 1, \dots, N \\ u(0; t_j) = k_1 & j = 1, \dots, M \\ u(x_N; t_j) = k_2 & j = 1, \dots, M \end{cases}$$

$$\begin{cases} S(x; t + \Delta t) = \varepsilon. [S(x_{i+1}, t_j) + S(x_{i-1}, t_j)] + (1 - 2\varepsilon) S(x_i, t_j) \\ \qquad + \dfrac{\varepsilon}{\alpha^2} \left( u(x_{i+1}, t_j) - 2u(x_i, t_j) + u(x_{i-1}, t_j) \right) \\ S(x_i; 0) = 0 & i = 1, \dots, N \\ S(0; t_j) = 0 & j = 1, \dots, M \\ S(x_N; t_j) = 0 & j = 1, \dots, M \end{cases}$$

donde $\varepsilon = \dfrac{\alpha^2 \Delta t}{\Delta x^2}$ es el parámetro de estabilidad numérica, el cual debe ser $\varepsilon < \dfrac{1}{2}$ para que resulte estable.

Las constantes elegidas para las simulaciones numéricas son:

Longitud $l$=0.4 m  
Fuente caliente $k_1$=100 °C  
Fuente fria $k_2$=0 °C  

*Condición inicial* Ta=25°C  
Deltax=0.01 m  
Deltat=0.1 Seg  

## **Resultados**

Con la finalidad de analizar la bondad del modelo se simularon los perfiles de temperatura para barras de distintos materiales. Se utilizó para la simulación una fuente caliente de 100°C en el extremo izquierdo y una fuente fría de 0°C en el extremo derecho. En todos los casos se consideró como condición inicial la temperatura ambiente (25°C). Las Figuras 1, 2 y 3 muestran el estado evolutivo en distintos puntos específicos de la barra (x = l/4, x = l/2 x = 3 l/4), mientras que la Figura 4 muestra la distribución espacial y temporal de temperaturas para un material particular (Cu).

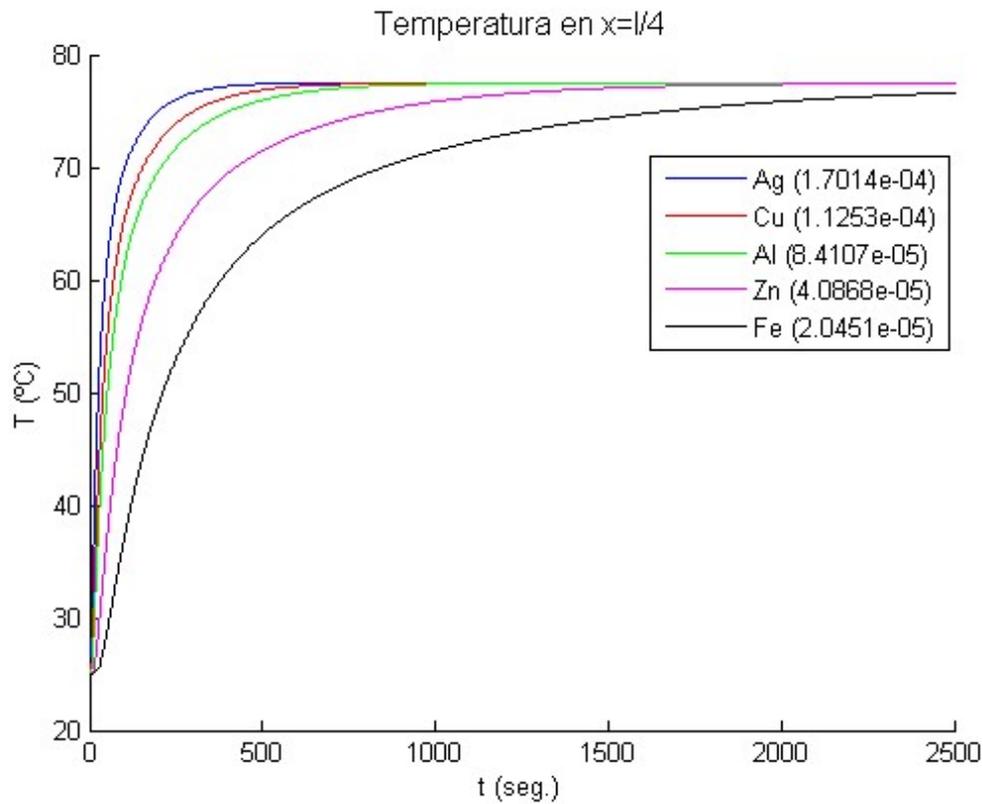

**Figura 1**: *Perfil temporal de temperatura para varios materiales en $x = l/4$.*

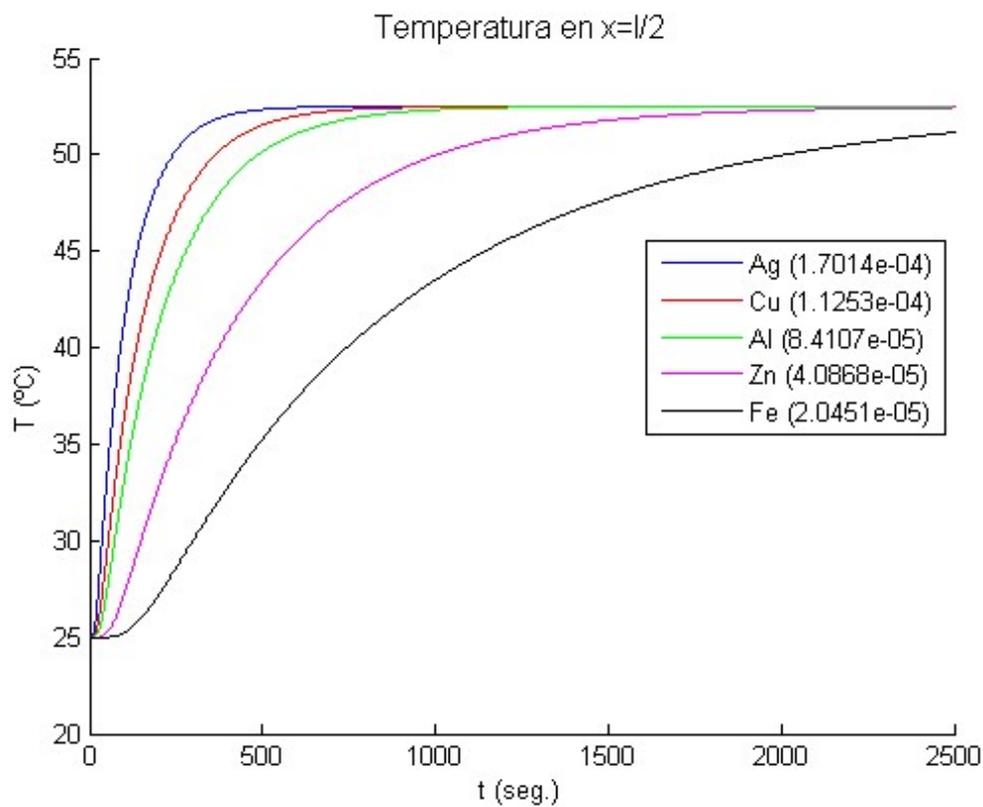

**Figura 2**: *Perfil temporal de temperatura para varios materiales en $x = l/2$.*

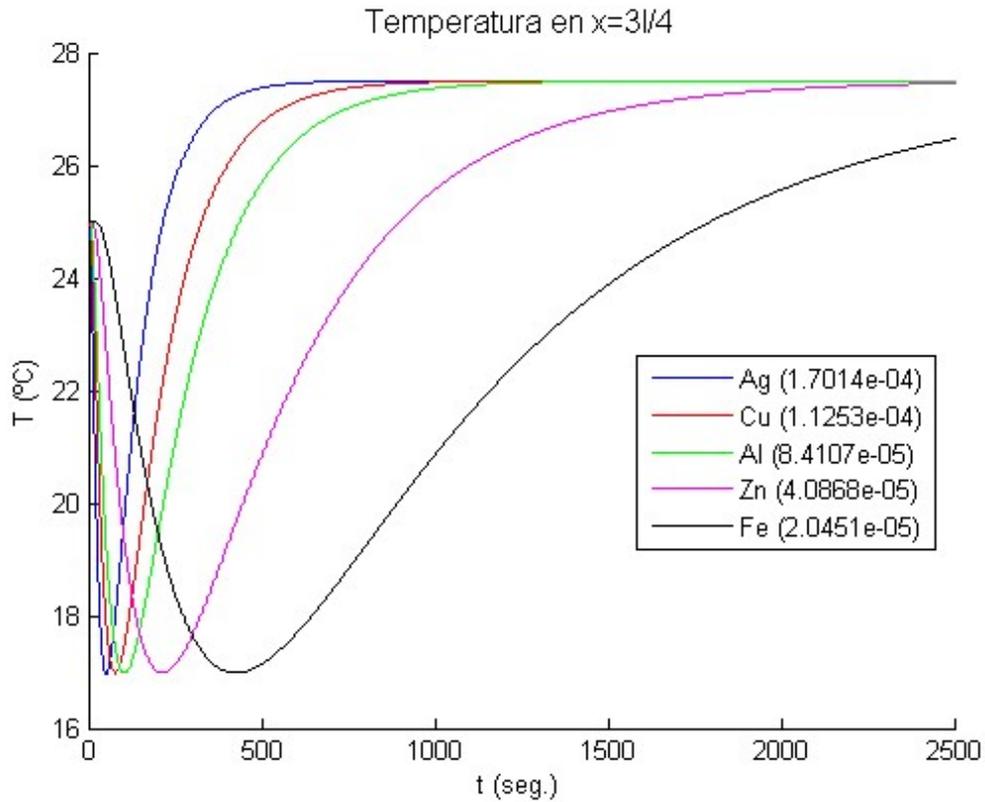

**Figura 3**: *Perfil temporal de temperatura para varios materiales en $x = 3\,l/4$.*

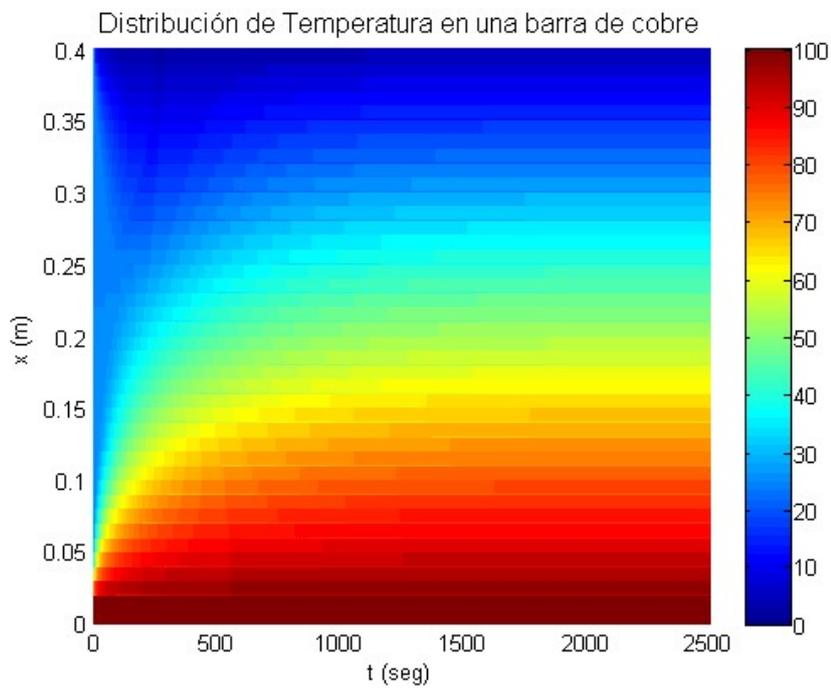

**Figura 4**: *Distribución de temperaturas en una barra de cobre.*

Por otra parte se graficó sensibilidad en función del tiempo para distintos materiales en las posiciones consideradas anteriormente. Estos gráficos se muestran en las Figuras 5-7.

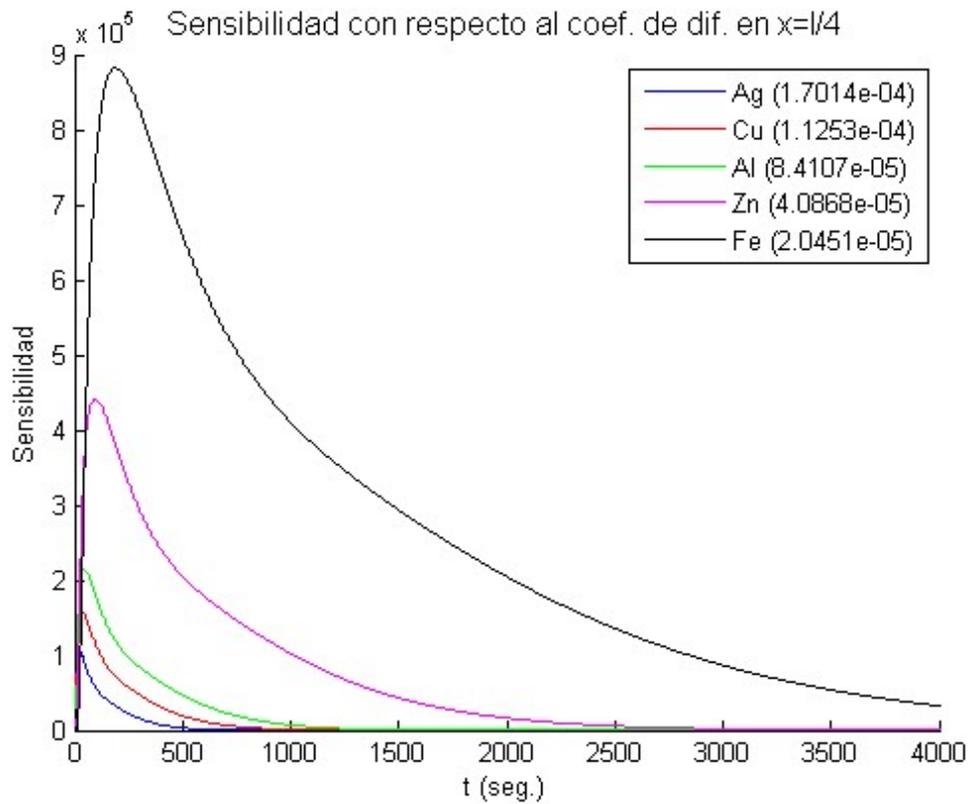

**Figura 5**: *Sensibilidad con respecto al coeficiente de difusividad en $x = l/4$.*

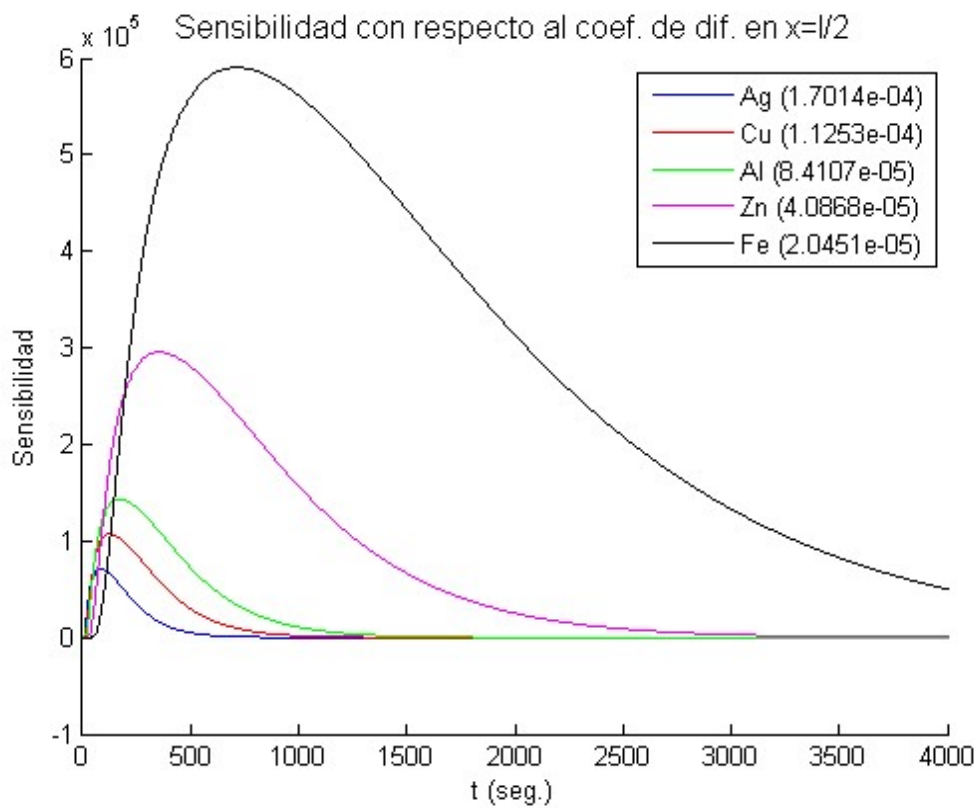

**Figura 6**: *Sensibilidad con respecto al coeficiente de difusividad en $x = l/2$.*

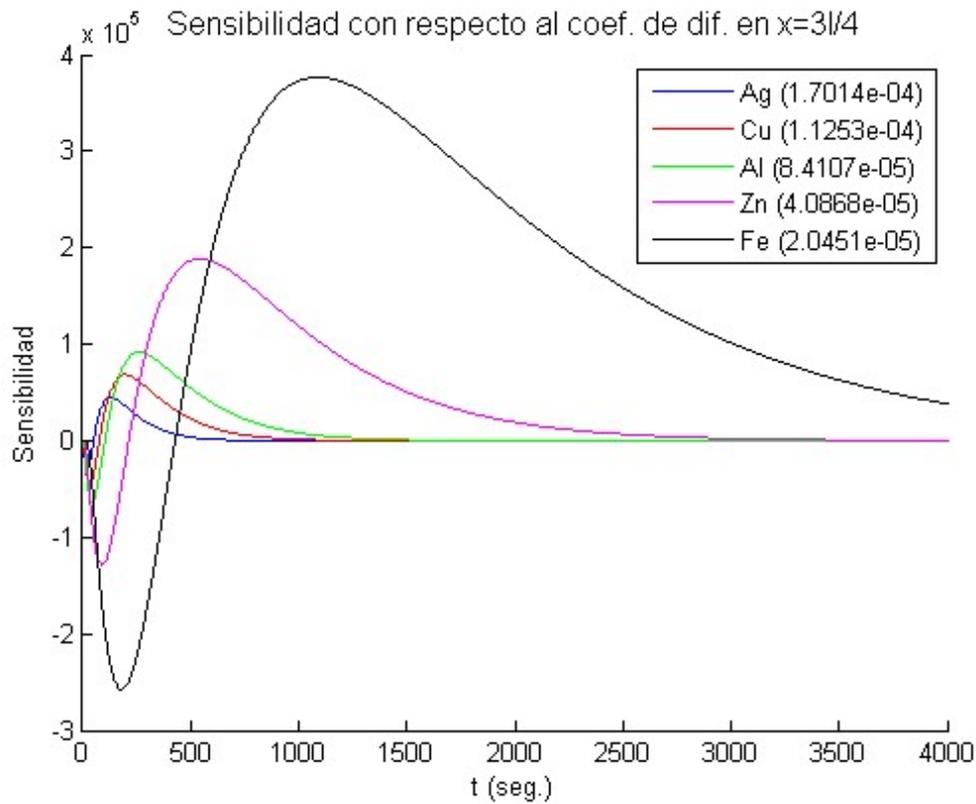

**Figura 7**: *Sensibilidad con respecto al coeficiente de difusividad en $x = 3\,l/4$.*

La distribución espacial y temporal de sensibilidad para el cobre (Cu) se puede apreciar en la Figura 8.

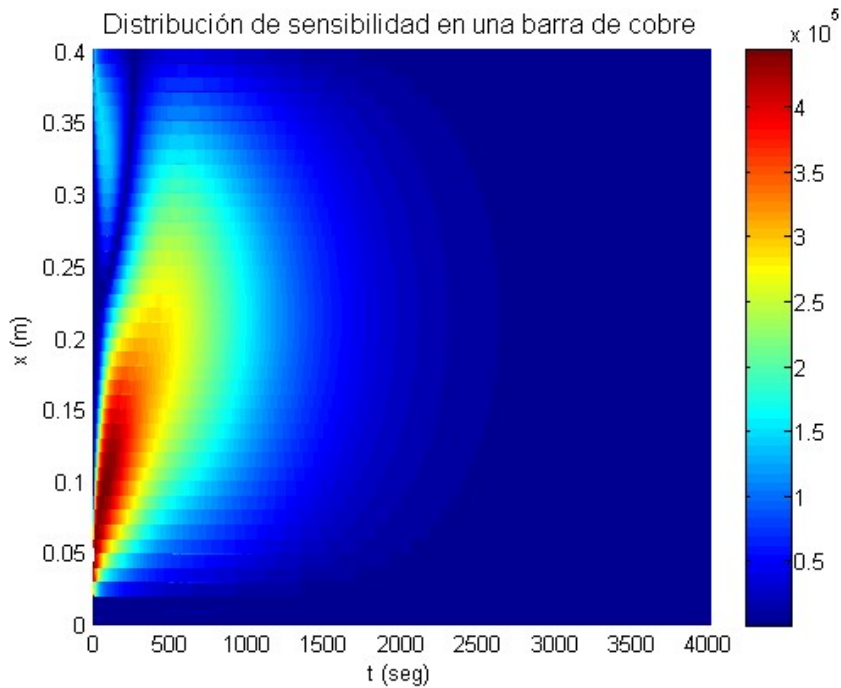

**Figura 8**: *Distribución de sensibilidad en una barra de cobre.*

Para la estimación del coeficiente de difusividad, consideramos $K=3$ mediciones de temperaturas en $x_0=l/2$ para distintos instantes de tiempo.

La estimación se hizo usando funciones de Matlab. La temperatura observada en cada caso fue simulada a partir de la solución numérica del problema directo a las cuales se las perturbó sumándole errores normalmente distribuidos con diferentes valores de desvío para simular errores de medición.

Se consideraron los valores correspondiente a diferentes metales: Plata, Aluminio, Zinc y Hierro.

Las Figuras 9-12 muestran los resultados obtenidos, los puntos de medición u observación de la temperatura.

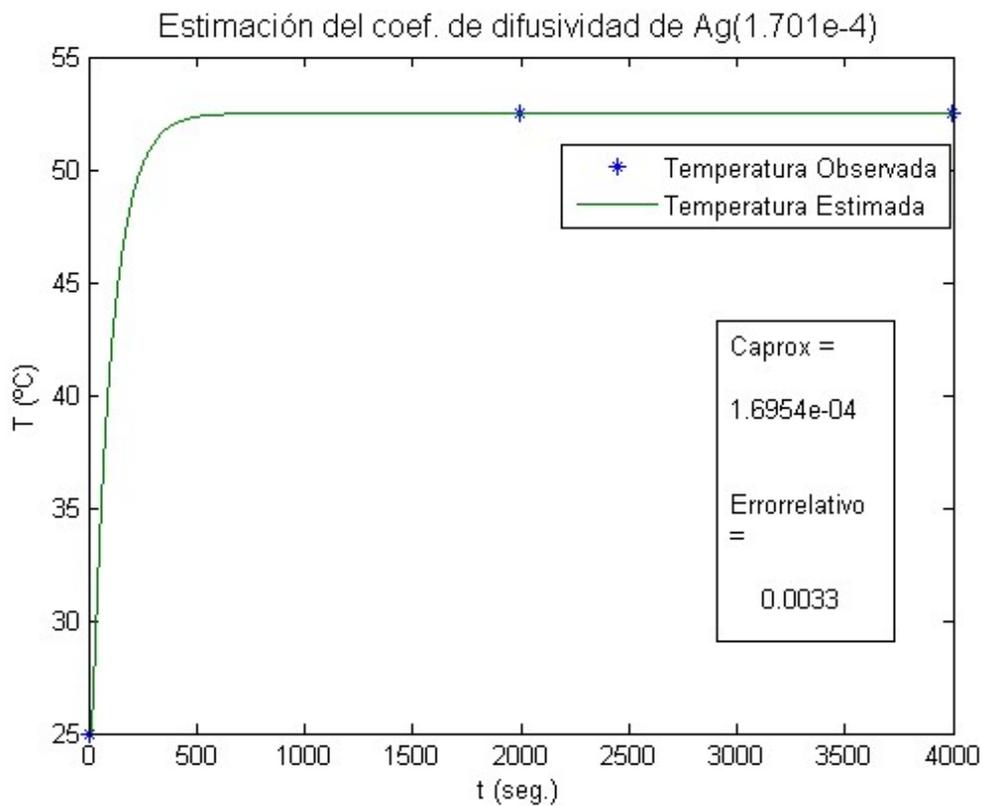

**Figura 9**: *Estimación del coeficiente de difusividad de la Plata.*

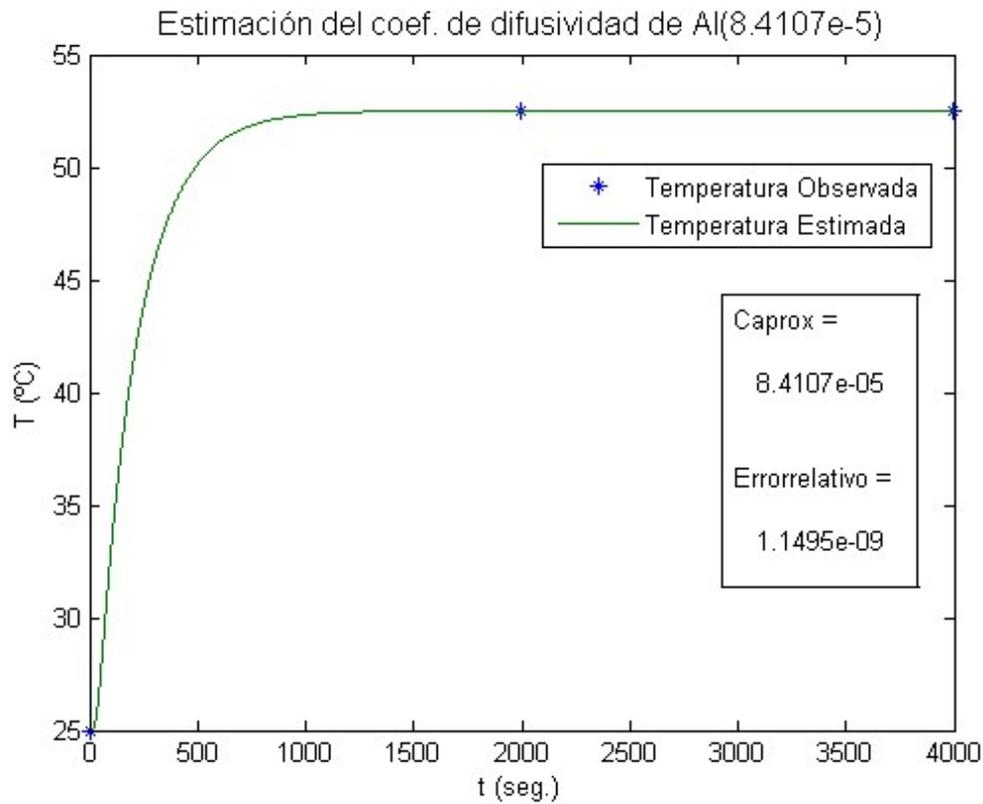

**Figura 10**: *Estimación del coeficiente de difusividad del Aluminio.*

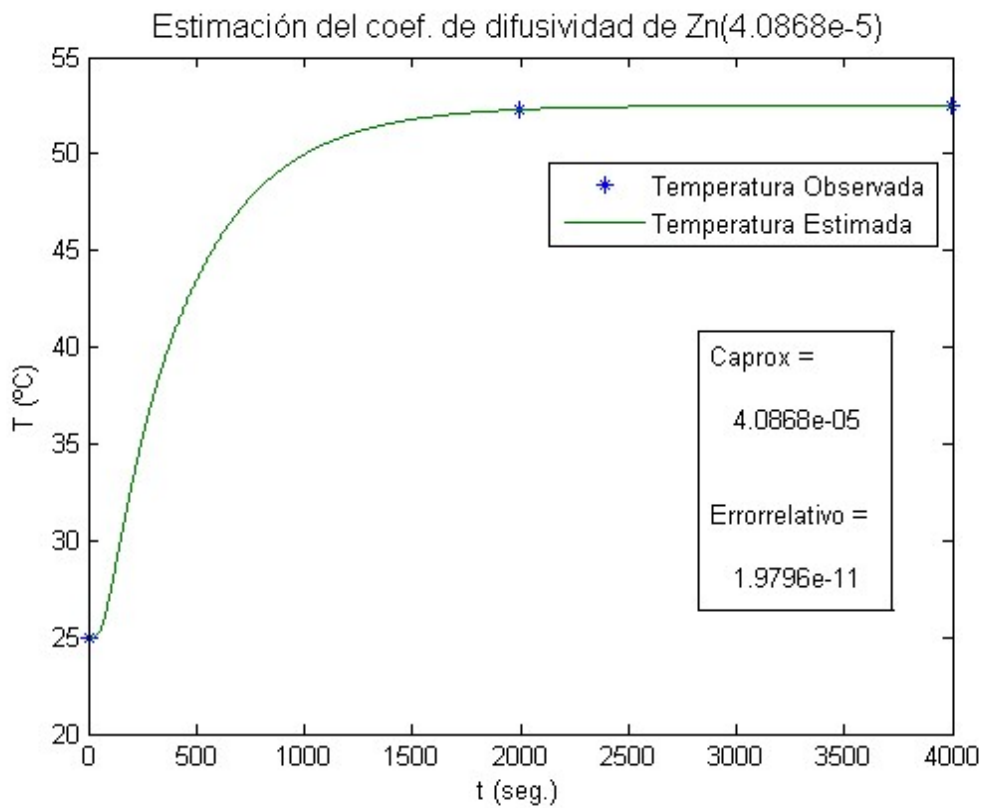

**Figura 12**: *Estimación del coeficiente de difusividad del Zinc.*

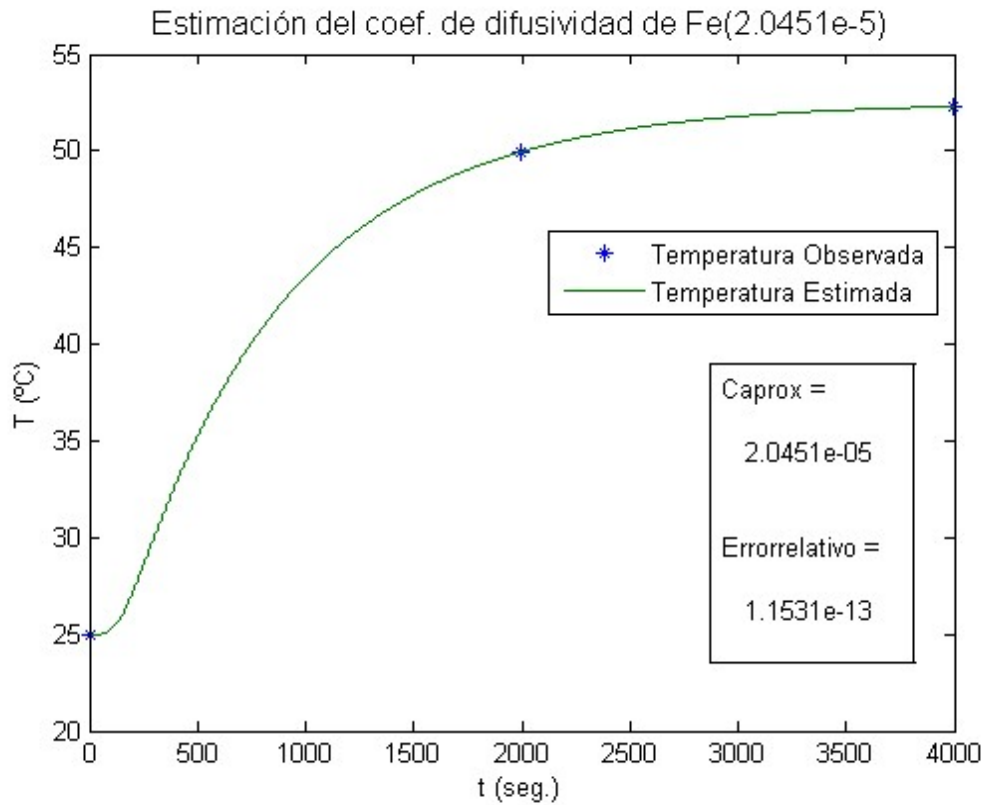

**Figura 12**: *Estimación del coeficiente de difusividad del Hierro.*

## Discusiones y conclusiones

En la Figura 1 se grafican los perfiles de temperatura en $x = l/4$ para 5 materiales diferentes. Podemos observar que en todos los casos la temperatura límite se aproxima a 78ºC (temperatura aproximada en el estado estacionario para $x = l/4$). Este valor límite es el mismo para todos los materiales debido a que en el estado estacionario todos tienen el mismo perfil de temperatura, ya que éste no depende del material sino de la posición con respecto a las fuentes y los valores de las mismas. Resultados análogos se observan para la Figura 2 y la Figura 3, siendo las temperaturas límite aproximadamente 53ºC (para $x = l/2$) y 27ºC (para $x = 3l/4$), respectivamente. Recordemos que en $x = 0$ se encuentra la fuente de 100°C mientras que en $x = l$ la temperatura es 0°C, esto explica la diferencia en los valores límite de temperatura. Más aún, en el estado estacionario, el perfil será una recta que une esos dos valores de temperatura, por lo que la temperatura en cada posición dependerá únicamente de la distancia a las fuentes.

Además, en la Figura 3, se observa un decrecimiento de la temperatura en los primeros instantes, esto se debe a la cercanía con la fuente fría (0°C), menor que la condición inicial (25°C). En la Figura 4 se muestran los valores que se obtienen en el caso del cobre en función del tiempo y de la posición.

Por otro lado se puede observar que los materiales más difusivos (coeficiente de difusividad más alto) alcanzan el estado estacionario más rápido, esto es consistente con los resultados teóricos.

Con respecto a la sensibilidad, la solución será más influenciada por el valor del coeficiente de difusividad, cuanto mayor es el valor de sensibilidad (en valor absoluto). Observando las Figuras 5-7, lo primero que se observa es que en todos los casos, los valores de sensibilidad decrecen con el tiempo, acercándose a 0. Por otro lado, la solución es más sensible a perturbaciones en el valor del parámetro para materiales menos difusivos. Los valores más altos de sensibilidad se encuentran en los primeros instantes de tiempo, y cercanos a la fuente caliente. En la Figura 8 se muestran los valores que se obtienen en el caso del cobre. La situación es similar para los diferentes materiales.

Se observa, por otro lado, que la sensibilidad de la temperatura con respecto al parámetro es nula cuando la temperatura alcanza el valor estacionario.

Por último, con respecto a la estimación del parámetro, se observa en todos los casos una buena precisión en la aproximación utilizando tan sólo 3 datos, en $t_1$=0 seg. , $t_2$=2000 seg., $t_3$=4000 seg, para $x_0$=l/2. Las Figuras 9-12 muestran las observaciones junto con la curva que se obtiene usando el valor estimado, además de los valores de las estimaciones y sus correspondientes errores relativos.